
\documentstyle{amsppt}
\magnification=\magstephalf
\font\medium=cmbx10 scaled \magstep1

\vsize=7.75truein
\hsize=6.25truein
\baselineskip=18pt
\parindent=20pt
\TagsOnRight

\def\IC{{\Bbb C}}

\def\IZ{{\Bbb Z}}
\def\IP{{\Bbb P}}

\def\spb{\smallpagebreak}
\def\mpb{\vskip 0.5truecm}
\def\bpb{\vskip 1truecm}
\def\added{{\it Added in second edition\/}:\ }
\vcorrection{-0.2 truecm}
\hcorrection{0.5 truecm}
\NoBlackBoxes
\def\folio{\tenrm\ifnum\pageno<0 \lowercase\expandafter{\romannumeral-\pageno}
 \else\number\pageno \fi}
\leftheadtext\nofrills{S. Katz}
\rightheadtext\nofrills{Rational Curves on Calabi--Yau Threefolds}



\topspace{3truecm}

\medium
\centerline{Rational Curves on Calabi--Yau Threefolds}

\vskip 2truecm

\rm
\centerline{Sheldon Katz\footnote{Supported in part by NSF grant DMS-9311386.}}

\spb
\centerline{Department of Mathematics}
\centerline{Oklahoma State University}
\centerline{Stillwater, OK 74078}

\bpb

In the conformal field theory arising from the compactification of
strings on a Calabi--Yau threefold $X$, there naturally arise fields
corresponding to harmonic forms of types $(2, 1)$ and $(1, 1)$ on $X$ [6].
The uncorrected Yukawa couplings on $H^{2, 1}$ and $H^{1, 1}$ are cubic
forms that can be constructed by techniques of algebraic geometry -- [31]
contains a nice survey of this in a general context in a language written
for mathematicians.  The cubic form on the space $H^{p, 1}$ of harmonic
$(p, 1)$ forms is given by the intersection product $(\omega_i, \omega_j,
\omega_k) \mapsto \int_X \omega_i \land \omega_j \land \omega_k$ for
$p = 1$, while for $p = 2$ there is a natural formulation in terms of
infinitesimal variation of Hodge structure [31, 9].  The Yukawa couplings
on $H^{2, 1}$ are exact, while those on $H^{1, 1}$ receive instanton
corrections.  In this context, the {\it instantons} are non-constant
holomorphic maps $f: \IC\IP^1 \to X$.  The image of such a map is a
{\it rational curve\/} on $X$, which may or may not be smooth.  If
the rational curve $C$ does not move inside $X$, then the contribution
of the instantons which are generically $1-1$ (i.e. birational) maps
with image $C$ can be written down explicitly -- this contribution only
depends on the topological type of $C$, or more or less equivalently,
on the integrals $\int_C f^* J$ for $(1, 1)$ forms $J$ on $X$.

\spb
If the conformal field theory could also be expressed in terms of a
``mirror manifold'' $X'$, then the uncorrected Yukawa couplings on
$H^{2, 1}(X')$ would be the same as the corrected Yukawa couplings
on $H^{1, 1} (X)$.  So if identifications could be made properly,
the infinitesimal variation of Hodge structure on $X'$ would give
information on the rational curves on $X$.  In a spectacular paper
[7], Candelas et.\ al.\ do this when $X$ is a quintic threefold.
Calculating the Yukawa coupling on $H^{2, 1} (X')$ as the complex
structure of $X'$ varies gives the Yukawa coupling on $H^{1, 1} (X)$
as $J$ varies.  The coefficients of the resulting Fourier series are
then directly related to the instanton corrections.\footnote{\added
There are now many more examples where this has been
carried out, e.g.\ [8,16,20,27,30,32].}

\spb
To the mathematician, there are some unanswered questions in deducing
the number of rational curves of degree $d$ on $X$ from this [31].  I
merely cite one problem here, which I will state a little differently.
It is not yet known how to carry out the calculation of the instanton
correction associated to a continuous family of rational curves on
$X$.  If there were continuous families of rational curves on a
general quintic threefold $X$, this would complicate the instanton
corrections to the Yukawa coupling on the space of $(1, 1)$ forms.
The question of whether or not such families exist has not yet been
resolved.  But even worse, this complication is present in any case --
any degree $m$ mapping of $\IC\IP^1$ to itself may be composed with any
birational mapping of $\IC\IP^1$ to $X$ to give a new map from $\IC\IP^1$
to $X$; and this family has more moduli than mere reparametrizations of
the sphere.  In [7] this was explicitly recognized; the assertion
was made that such a family counts $1/m^3$ times.  This has since been
verified by Aspinwall and Morrison [3].\footnote{\added
Y.~Ruan has an interesting approach to giving a mathematically
precise meaning to this statement [35].}

\spb
My motivation in writing this note is to give a general feel for the
mathematical meaning of ``the number of rational curves on a Calabi--Yau
threefold'', in particular, how to ``count'' a family of rational curves.
I expect that these notions will directly correspond to the not as yet
worked out procedure for calculating instanton corrections associated to
general families of rational curves.  More satisfactory mathematical
formulations are the subject of work still in progress.

\spb
I want to take this opportunity to thank David Morrison for his suggestions,
to the Mathematics Department at Duke University for its hospitality
while this manuscript was being prepared, and to the organizers of the
Mirror Symmetry Workshop for providing a wonderful opportunity for
mathematicians and physicists to learn from each other.

\bpb
\subhead\nofrills 1. Rational Curves, Normal Bundles, Deformations
\endsubhead

\mpb
Consider a Calabi--Yau threefold $X$ containing a smooth rational curve
$C \simeq \IC\IP^1$.  The normal bundle $N_{C/X}$ of $C$ in $X$ is
defined by the exact sequence
$$
0 \to T_C \to T_X|_C \to N_{C/X} \to 0 \ . \tag1
$$
%
$N = N_{C/X}$ is a rank $2$ vector bundle on $C$, so $N = \Cal O(a)
\oplus \Cal O(b)$ for some integers $a, b$.  Now $c_1(T_C) = 2$,
and $c_1(T_X) = 0$ by the Calabi--Yau condition.  So the exact sequence
(1) yields $c_1(N) = -2$, or $a + b = -2$.  One ``expects'' $a = b = -1$
in the general case.  This is because there is a moduli space of
deformations of the vector bundle $\Cal O(a) \oplus \Cal O(b)$, and the
general point of this moduli space is $\Cal O(-1) \oplus \Cal O(-1)$,
no matter what $a$ and $b$ are, as long as $a + b = -2$.

\spb
Let $\Cal M$ be the moduli space of rational curves in $X$.  The tangent
space to $\Cal M$ at $C$ is given by $H^0(N)$ [28, \S12].  In other
words, $\Cal M$ may be locally defined by finitely many equations in
$\dim H^0(N)$ variables.

\spb
\proclaim{Definition} $C$ is {\it infinitesimally rigid\/} if $H^0(N) = 0$.
\endproclaim

\spb
Infinitesimal rigidity means that $C$ does not deform inside $X$, not
even to first order.  Note that $H^0(N) = 0$ if and only if
$a = b = -1$.  Thus

\item{$\bullet$} $C$ is infinitesimally rigid if and only if $a = b = -1$.
\item{$\bullet$} $C$ deforms, at least infinitesimally, if and only if
$(a, b) \ne (-1, -1)$.

\spb
$\Cal M$ can split up into countably many irreducible components.
For instance, curves with distinct homology classes in $X$ will lie
in different components of $\Cal M$.  However, there will be at most
finitely many components of $\Cal M$ corresponding to rational curves
in a fixed homology class.

\spb
There certainly exist Calabi--Yau threefolds $X$ containing positive
dimensional families of rational curves.  For instance, the Fermat
quintic threefold $x^5_0 + \cdots + x^5_4 = 0$ contains the family
of lines given parametrically in the homogeneous coordinates $(u, v)$
of $\IP^1$ by $(u, -u, av, bv, cv)$, where $(a, b, c)$ are the parameters
of the plane curve $a^5 + b^5 + c^5 = 0$.  However, suppose that all
rational curves on $X$ have normal bundle $\Cal O(-1) \oplus \Cal O(-1)$.
Then since $\Cal M$ consists entirely of discrete points, the remarks
above show that there would be only finitely many rational curves in $X$
in a fixed homology class.  Since a quintic threefold has $H_2(X, \IZ)
\simeq \IZ$, the degree of a curve is essentially the same as its
homology class.  This discussion leads to Clemens' conjecture [11]:

\spb
\proclaim{Conjecture (Clemens)} A general quintic threefold contains
only finitely many rational curves of degree $d$, for any $d \in \IZ$.
These curves are all infinitesimally rigid.
\endproclaim

\spb
Clemens' original constant count went as follows [10]: a rational curve of
degree $d$ in $\IP^4$ is given parametrically by $5$ forms $\alpha_0(u, v),
\ldots, \alpha_4(u, v)$, each homogeneous of degree $d$ in the homogeneous
coordinates $(u:v)$ of $\IP^1$.  These $\alpha_i$ depend on $5(d+1)$
parameters.  On the other hand, a quintic equation $F(x_0, \ldots, x_4) = 0$
imposes the condition $F\big(\alpha_0(u, v), \ldots, \alpha_4(u, v)\big)
\equiv 0$ for the parametric curve to be contained in this quintic
threefold.  This is a polynomial equation of degree $5d$ in $u$ and $v$.
Since a general degree $5d$ polynomial $\Sigma a_i u^i v^{5d-i}$ has
$5d+1$ coefficients, setting these equal to zero results in $5d+1$
equations among the $5(d+1)$ parameters of the $\alpha_i$.  If $F$ is
{\it general\/}, it seems plausible that these equations should impose
independent conditions, so that the solutions should depend on $5d+5-(5d+1)
= 4$ parameters.  However, any curve has a $4$-parameter family
of reparametrizations $(u, v) \mapsto (au + bv, cu + dv)$, so there are
actually a zero dimensional, or finite number, of curves on the general
$F = 0$.

\spb
The conjecture is known to be true for $d \le 7$ [23].\footnote{\added
This result has recently been extended to $d=8$ and $d=9$ by T.~Johnsen
and S.~Kleiman [21].}
  For any $d$, it can
even be proven that there exists an infinitesimally rigid curve of degree
$d$ on a general $X$.  Similar conjectures can be stated for other
Calabi--Yau threefolds.\footnote{\added The existence of an infinitesimally
rigid
curve of arbitrary degree $d$ on a general $(2,4)$ complete intersection
Calabi-Yau threefold in $\IP^5$ has been proven recently by K.~Oguiso [33].}

\spb
There are many kinds of non-rational curves which appear to occur in
finite number on a general quintic threefold.  For instance, elliptic
cubic curves are all planar.  The plane $P$ that one spans meets the
quintic in a quintic curve containing the cubic curve.  The other
component must be a conic curve. This sets up a $1-1$ correspondence
between elliptic cubics and conics on any quintic threefold.  Hence the
number of elliptic cubics on a general quintic must be the same
as the number of conics, $609250$ [23].  Finiteness of elliptic quartic
curves has been proven by Vainsencher [37]; the actual number has not
yet been computed.\footnote{\added This
number has recently been computed by G.~Ellingsrud and S.A.~Str\o mme
to be 3718024750 (private communication).  They use the description
of the parameter space for elliptic curves given by Avritzer-Vainsencher
and a formula of Bott, whose use they explained in [15].
Also, a method for predicting instantons of
higher genus has been introduced by Bershadsky, Cecotti, Ooguri, and Vafa
in a series of two papers [4,5].  The predictions are compatible with the
calculation of Ellingsrud and Str\o mme as well as with previously known
numbers.  This method has been recently extended to multiparameter families
as part of a paper by Candelas, de la Ossa, Font, Katz, and
Morrison; the method has been applied again in a companion paper [8].}

\spb
On the other hand, there are infinitely many plane quartics on {\it any\/}
quintic threefold: take any line in the quintic, and each of the
infinitely many planes containing the line must meet the quintic in the
original line union a quartic.

\spb
If a curve has $N \simeq \Cal O \oplus \Cal O(-2)$, then $C \subset X$
deforms to first order.  In fact, since $H^0(N)$ is one-dimensional,
there is a family of curves on $X$ parametrized by a single variable
$t$, subject to the constraint $t^2 = 0$.  In other words, start with a
rational curve given parametrically by forms $\alpha_0, \ldots, \alpha_4$,
homogeneous of degree $d$ in $u$ and $v$.  Take a perturbation $\alpha_i
(u, v; t) = \alpha_i(u, v) + t \alpha_i'(u, v)$, still homogeneous in
$(u, v)$.  Form the equation $F(\alpha_0, \ldots, \alpha_4) = 0$ and
formally set $t^2 = 0$; the resulting equation has a $5$ dimensional
space of solutions for the $\alpha'_i$, which translates into a
unique solution up to multiples and reparametrizations of $\IP^1$.
The curve $C$ but may or may not deform to second order.  $C$ deforms
to $n^{\text{th}}$ order for all $n$ if and only if $C$ moves in a
$1$-parameter family.  A pretty description of the general situation
is given in [34].

\spb
If $C$ deforms to $n^{\text{th}}$ order, but not to $(n + 1)^{\text{th}}$
order, then one sees that while $C$ is an isolated point in the
moduli space of curves on $X$, it more naturally is viewed
as the solution to the equation $t^{n+1} = 0$ in one variable $t$.  So
$C$ should be viewed as a rational curve on $X$ with multiplicity
$n+1$.

\spb
If a curve has $N \simeq \Cal O(1) \oplus \Cal O(-3)$, $C$ has a $2$
parameter space of infinitesimal deformations, and the structure of
$\Cal M$ at $C$ is correspondingly more complicated.  An example is given
in the next section.  The general situation has not yet been worked
out.

\bpb

\subhead\nofrills 2. Counting Rational Curves
\endsubhead

\mpb
In this section, a general procedure for calculating the number of
smooth rational curves of a given type is described.  Alternatively,
a canonical definition of this number can be given using the Hilbert
scheme (this is what was done by Ellingsrud and Str\o mme in their
work on twisted cubics [14]); however, it is usually quite difficult
to implement a calculation along these lines.

\spb
Embed the Calabi--Yau threefold $X$ in a larger compact space
$\IP$ (which may be thought of as a projective space, a weighted
projective space, or a product of such spaces).  $\Cal M_\lambda$
will denote the moduli space parametrizing smooth rational curves
in $\IP$ of given topological type or degree $\lambda$.  $\hbox{Def}(X)$
denotes the irreducible component of $X$ in the moduli space of
Calabi--Yau manifolds in $\IP$ (here the K\"ahler structure is ignored).
In other words, $\hbox{Def}(X)$ parametrizes the deformations of $X$ in
$\IP$.

\spb
\item{(1)} Find a compact moduli space $\overline  {\Cal M}_\lambda$ containing
$\Cal  M_\lambda$ as a dense open subset, such that the points of
$\overline  {\Cal M}_\lambda - \Cal M_\lambda$ correspond to degenerate curves
of type $\lambda$.  $\overline  {\Cal M}_\lambda$ parametrizes degenerate
deformations of the smooth curve (not the mapping from $\IC\IP^1$
to $X$).  It is better for $\overline  {\Cal M}_\lambda$ to be smooth.

\item{(2)} Find a rank $r = \dim (\Cal M_\lambda)$ vector bundle $\Cal B$
on $\overline  {\Cal M}_\lambda$ such that

\itemitem{(a)} To each $X' \in \hbox{Def}(X)$ there is a section $s_{X'}$
of $\Cal B$ which vanishes at $C \in \overline  {\Cal M}_\lambda$ if and only
if $C \subset X'$.

\itemitem{(b)} There exists an $X' \in \hbox{Def}(X)$ such that $s_{X'}(C) = 0$
if and only if $C \in \Cal M_\lambda$ and $C \subset X'$.

\itemitem{(c)} $C$ is an isolated zero of $s_{X'}$ with $\text{mult}_C
(s_{X'}) = 1$ if and only if $N_{C/X} \simeq \Cal O(-1) \oplus \Cal O(-1)$.

\spb
\proclaim{Working Definition} The number of smooth rational curves
$n_\lambda$ of type $\lambda$ is given by the $r^{\text{th}}$ Chern
class $c_r(\Cal B)$.
\endproclaim

\spb
Why is this reasonable definition?  Suppose that an $X' \in \hbox{Def}(X)$
can be found with the properties required above, with the additional
property that there are only finitely many curves of type $\lambda$ on
$X'$, and that they are all infinitesimally rigid.  Then it can be
checked that the number of (possibly degenerate) curves of type $\lambda$
on $X'$ is independent of the choice of $X'$ satisfying the above
properties, and is also equal to $c_r(\Cal B)$.  This last follows since
$c_{\text{rank}(E)}(E)$ always gives the homology class of the zero locus
$Z$ of any section of any bundle $E$ on any variety $Y$, whenever
$\dim (Z) = \dim (Y) - \text{rank} (E)$.  In our case, $0 = \dim$
(\{{\it lines\/}\}) = $\dim (\Cal M_\lambda) - \text{rank} (\Cal B)$.
In essentially all known cases, the number of curves has been worked out
by the method of this working definition.  Examples are given below.

\spb
There is a potential problem with this working definition.  For families
of Calabi--Yau threefolds such that no threefold in the family contains
finitely many curves of given type, it may be that the ``definition''
depends on the choice of compactification and/or vector bundle, i.e.
this is not well-defined.  My reason for almost calling this method
a definition is that it {\it does\/} give a finite number corresponding
to an infinite family of curves, which {\it is\/} well-defined in the
case that the Calabi--Yau threefold in question belongs to a family
containing some other Calabi--Yau threefold with only finitely many
rational curves of the type under consideration.

\spb
In the case where the general $X$ contains irreducible singular
curves which are the images of maps from $\IP^1$, a separate but similar
procedure must be implemented to calculate these, since they give
rise to instanton corrections as well.  For example, there is a
$6$-parameter family of two-planes in $\IP^4$.  A two plane $P$ meets a
quintic threefold $X$ in a plane quintic curve.  For general $P$ and
$X$, this curve is a smooth genus $6$ curve.  But if the curve
acquired $6$ nodes, the curve would be rational.  This being $6$
conditions on a $6$ parameter family, one expects that a general $X$
would contain finitely many $6$-nodal rational plane quintic curves,
and it can be verified that this is indeed the case.\footnote{\added The number
of such curves has recently been calculated to be 17,601,000 by
I.~Vainsencher [38].}
The problem
is that there is no way to deform a smooth rational curve to such
a singular curve -- the dimension of $H^1({\Cal O}_C)$ is zero for a
smooth rational curve $C$, but positive for such singular curves,
and this dimension is a deformation invariant [19, Theorem 9.9].
For example, if one tried to deform a smooth twisted cubic curve to a
singular cubic plane curve by projecting onto a plane, there would
result an ``embedded point'' at the singularity, creating a sort of
discontinuity in the deformation process [19, Ex. 9.8.4].

\spb
\noindent
{\bf Examples:}

\mpb
\item{1.} Let $X \subset \IP^4$ be a quintic threefold.  Take $\lambda = 1$,
so that we are counting lines.  Here $\Cal M_1 = G(1, 4)$ is the
Grassmannian of lines in $\IP^4$ and is already compact,  so take
$\overline  {\Cal M}_1 = \Cal M_1$.  Let $\Cal B = \text{Sym}^5 (U^*)$, where
$U$, the universal bundle, is the rank $2$ bundle on $\Cal M_1$ whose
fiber over a line $L$ is the $2$-dimensional subspace $V \subset \IC^5$
yielding $L \subset \IP^4$ after projectivization.  Note that
$\text{rank}(\Cal B) = \dim (\Cal M_1) = 6$.  $\hbox{Def}(X) \subset \IP
(H^0(\Cal O_{\IP^4}(5)))$ is the subset of smooth quintics.  A quintic
$X$ induces a section $s_X$ of $\Cal B$, since an equation for $X$
is a quintic form on $\IC^5$, hence induces a quintic form on $V$ for
$V \subset \IC^5$ corresponding to $L$.  Clearly $s_X(L)=0$ if and only
if $L \subset X$.  The above conditions are easily seen to hold.
$c_6(\Cal B) = 2875$ is the number of lines on $X$.  This calculation
is essentially the same as that done for cubics in [2, Thm.\ 1.3], where
dual notation is used, so that the $U^*$ used here becomes the universal
quotient bundle $Q$ in [2].  The number $2875$ agrees with the result
of Candelas et.\ al.\ [7].

\item{(2)} Continuing with the quintic, take $\lambda = 2$.  Any conic
$C$ necessarily spans a unique $2$-plane containing $C$.  Let $G = G(2, 4)$
be the Grassmannian of $2$-planes in $\IP^4$, and let $U$ be the rank $3$
universal bundle on $G$.  Put $\overline  {\Cal M}_2 = \IP(\text{Sym}^2(U^*))$
be the projective bundle over $G$ whose fiber over a plane $P$ is the
projective space of conics in $P$.  Clearly $\Cal M_2 \subset \overline {\Cal
M}_2$
(but they are not equal; $\overline  {\Cal M}_2$ also contains the union
of any two lines or a double line in any plane).  Let $\Cal B =
\text{Sym}^5(U^*)/(\text{Sym}^3(U^*) \otimes \Cal O_{\IP}(-1))$ be the
bundle on $\overline  {\Cal M}_2$ of quintic forms on the $3$ dimensional
vector space $V \subset \IC^5$, modulo those which factor as any cubic
times the given conic.  $(\Cal O_{\IP}(-1)$ is the line bundle whose
fiber over a conic is the one dimensional vector space of equations for
the conic within its supporting plane.  The quotient is relative to the
natural embedding $\text{Sym}^3(U^*) \otimes \Cal O_{\IP}(-1) \to
\text{Sym}^5(U^*)$ induced by multiplication.)  $\text{rank} (\Cal B) =
\dim (\overline  {\Cal M}_2) = 11$.  $c_{11}(\Cal B) = 609250$.  See [23] for
more
details, or [12] for an analogous computation in the case of quartics.
The number $609250$ agrees with the result of Candelas et.\ al.

\item{(3)} Again consider the quintic, this time with $l = 3$.  In [14],
Ellingsrud and Str\o mme take $\overline  {\Cal M}_3$ to be the closure
of the locus of smooth twisted cubics in $\IP^4$ inside the Hilbert
scheme.  This space has dimension $16$.  $\Cal B$ is essentially the
$\text{rank}\ 16$ bundle of degree $15$ forms on $\IP^1$ induced from
quintic polynomials in $P^4$ by the degree $3$ parametrization of
the cubic (it must be shown that this makes sense for degenerate
twisted cubics as well).  The equation of a general quintic gives a
section of $\Cal B$, vanishing precisely on the set of cubics
contained in $X$.  $c_{16} (\Cal B) = 317206375$, again agreeing
with the result of Candelas et.\ al.

\spb
The key to the first two calculations are the Schubert calculus for
calculating in Grassmannians [18, Ch.\ 1.5] and standard formulas for
projective bundles [19, Appendix A.3].  The third calculation is more
intricate.

\spb
Regarding complete intersection Calabi--Yau manifolds, similar examples
are found in [29,24] for lines and [36] for conics.\footnote{\added
These numbers agree with the numbers calculated by A.~Libgober and
J.~Teitelbaum using proposed mirrors for these manifolds.}

\spb
Note that the number of curves $c_r (\Cal B)$ in no way depends on
the choice of $X$, even if $X$ is a degenerate Calabi--Yau threefold,
or contains infinitely many rational curves of type $\lambda$.  It turns
out that a natural meaning can be assigned to this number.\footnote{\added
Y.~Ruan has shown along these lines that this number can be
given a precise meaning by deforming $X$ to an almost complex manifold [35].}
In fact,
the moduli space of curves of type $\lambda$ on $X$ splits up into
``distinguished varieties'' $Z_i$ [17], and a number, the
{\it equivalence\/} of $Z_i$, can be assigned to each distinguished
variety (the number is $1$ for each $\Cal O(-1) \oplus \Cal O(-1)$
curve).  This number is precisely equal to the number of curves on
$X$ in $Z_i$ which arise as limits of curves on $X'$ and $X'$ approaches
$X$ in a $1$-parameter family [17, Ch.\ 11].

\spb
\noindent
{\bf Examples:}

\mpb
\item{(1)} If a quintic threefold is a union of a hyperplane and a quartic,
then the quintic contains infinitely many lines and conics.  However,
in the case of lines, given a general $1$-parameter family of quintics
approaching this reducible quintic, $1275$ lines approach the hyperplane,
and $1600$ lines approach the quartic [22].  So the infinite set of lines
in the hyperplane ``count'' as $1275$, while those in the quartic count
as $1600$.  For counting conics, the $609250$ conics distribute
themselves as $187850$ corresponding to the component of conics in the
hyperplane, $258200$ corresponding to the component of conics in the
quartic, and $163200$ corresponding to the component of conics which
degenerate into a line in the hyperplane union an intersecting line in
the quartic [25].  Note that this reducible conic lies in $\overline  {\Cal
M}_2
- \Cal M_2$, and illustrates why $\Cal M_2$ itself is insufficient
for calculating numbers when there are infinitely many curves.  Most
of these numbers have been calculated recently by Xian Wu [39] using a
different method.

\item{(2)} If a quintic threefold is a union of a quadric and a cubic,
the lines on the quadric count as $1300$, and the lines on the cubic
count as $1575$ [22, 39].  The conics on the quadric count as $215950$,
while the conics which degenerate into a line in the quadric union an
intersecting line in the cubic count as $609250-(215950+243900) = 149400$,
but this has not been checked directly yet.

\item{(3)} There are infinitely many lines on the Fermat quintic threefold
$x^5_0 + \cdots + x^5_4 = 0$.  These divide up into $50$ cones, a typical
one being the family of lines given parametrically in the homogeneous
coordinates $(u, v)$ of $\IP^1$ by $(u, -u, av, bv, cv)$, where $(a, b, c)$
satisfy $a^5 + b^5 + c^5 = 0$.  Each of these count as $20$.  There are
also $375$ special lines, a typical one being given by the  equations
$x_0 + x_1 = x_2 + x_3 = x_4 = 0$ (these lines were also noticed in [13]).
These lines $L$ count with multiplicity $5$.  Note that $50 \cdot 20 +
375 \cdot 5 = 2875$ [1].  This example illustrates the potential
complexity in calculating the distinguished varieties $Z_i$ -- some
components can be embedded inside others.  This may be understood as well
by looking at the moduli space of lines on the Fermat quintic locally
at a line corresponding to a special line.  Since $N_{L/X} \simeq \Cal O(1)
\oplus \Cal O(-3)$, $H^0(N_{L/X})$ has dimension $2$, and the moduli
space of lines is a subset of a $2$-dimensional space.  A calculation
shows that it is locally defined inside $2$ dimensional $(x, y)$ space
by the equations $x^2 y^3 = x^3 y^2 = 0$.  The $x$ and $y$ axes
correspond to lines on each of the $2$ cones, each occurring with
multiplicity $2$; there would be just one  equation $x^2 y^2 = 0$ if
the special line corresponding to $(0, 0)$ played no role; since this
is not the case, it can be expected to have its own contribution;
i.e. each special line is a distinguished variety.

\item{(4)} Examples for lines in complete intersection Calabi--Yau
threefolds were worked out in [24].

\spb
These numbers can also be calculated by intersection-theoretic techniques.
Let $s(Z_i, \overline  {\Cal M}_\lambda)$ denote the {\it Segre class\/}
of $Z_i$ in $\overline  {\Cal M}_\lambda$.  If $Z_i$ is smooth, this is
simply the formal inverse of the total Chern class $1 + c_1(N) + c_2(N) +
\ldots$ of the normal bundle $N$ of $Z_i$ in $\overline  {\Cal M}_\lambda$.
Then if $Z_i$ is a connected component of the zero locus of $s_X$, the
equivalence of $Z_i$ is the zero dimensional part of $c(\Cal B) \cap
s(Z_i, \overline  {\Cal M}_\lambda)$ [17, Prop. 9.1.1].

\spb
If $Z_i$ is an irreducible component which is not a connected component,
then this formula is no longer applicable.  However, I have had recent
success with a new method that supplies ``correction terms'' to this
formula.  The method is currently ad hoc (the most relevant success I
have had is in calculating the ``number'' of lines on a cubic surface
which is a union of three planes), but I expect that a more systematic
procedure can be developed.

\spb
This of course is a reflection of the situation in calculating instanton
corrections to the Yukawa couplings.  If there is a continuous family
of instantons, then calculating the corrections will be more difficult.
If the parameter space for instantons is smooth, this should make the
difficulties more manageable.  If the space is singular,
the calculation is more difficult.  I expect that the calculation
would be even more difficult if instantons occur in at least $2$
families that intersect.

\spb
Of course, in the calculation of the Yukawa couplings via path
integrals, there is no mention of vector bundles on the moduli space.
This indicates to me that there should be a mathematical definition
of the equivalence of a distinguished variety that does not refer to
an auxiliary bundle.
\footnote{\added
The current situation is that a good definition exists
as long as the dimension of $H^1(N)$ is independent
of all curves in the connected component [26].}
Along these lines, one theorem will be stated
without proof.

\spb
Let $Z$ be a $k$-dimensional unobstructed family of rational curves
on a Calabi--Yau threefold $X$.  There is the total space $\Cal Z
\subset Z \times X$ of the family, with projection map $\pi: \Cal Z
\to Z$ such that $\pi^{-1}(z)$ is the curve in $X$ corresponding to
$z$, for each $z \in Z$.  Let $N$ be the normal bundle of $\Cal Z$ in
$Z \times X$.  Define the equivalence $e(Z)$ of $Z$ to be the number
$c_k(R^1 \pi_* N)$.  For example, if $k = 0$, then the curve is
infinitesimally rigid, and $e(Z) = 1$.

\spb
\proclaim{Theorem}\footnote{\added
An explanation of this theorem appears in the second part of
[8].}
 Let $Z$ be an unobstructed family of rational
curves of type $\lambda$  on a Calabi--Yau threefold $X$.  Suppose
that $X$ deforms to a Calabi--Yau threefold containing only finitely
many curves of type $\lambda$.  Then precisely $e(Z)$ of these curves
(including multiplicity) approach curves of $Z$ as the Calabi--Yau
deforms to $X$.
\endproclaim

\spb
Some of the examples given earlier in this section can be redone via
this theorem.  Also, as anticipated by [7], it can be calculated that
a factor of $1/m^3$ is introduced by degree $m$ covers by a calculation
similar in spirit to that found in [3].

\newpage

\centerline{\bf References}

\mpb
\item{1.} A. Albano and S. Katz,  {\it Lines on the Fermat quintic threefold
and the infinitesimal generalized Hodge conjecture\/},  Trans.\ AMS {\bf 324}
(1991) 353--368.

\item{2.} A. Altman and S. Kleiman, {\it Foundations of the theory of Fano
schemes\/},  Comp.\ Math.\ {\bf 34} (1977) 3--47.

\item{3.} P. Aspinwall and D. Morrison,  {\it Topological field theory and
rational curves\/},  Comm.\ Math.\ Phys.\ {\bf 151} (1993) 245--262.

\item{4.} M.~Bershadsky, S.~Cecotti, H.~Ooguri and C.~Vafa.
{\it Holomorphic Anomalies
in Topological Field Theories\/}, with an appendix by S. Katz.
To appear in Nuc.\ Phys.\ B.

\item{5.} M.~Bershadsky, S.~Cecotti, H.~Ooguri and C.~Vafa.
{\it Kodaira-Spencer theory of gravity\/},
preprint, HUTP--93/A025.

\item{6.} P. Candelas, G. Horowitz, A. Strominger and E. Witten,
{\it Vacuum configurations for super-strings\/},  Nuc.\ Phys.\
{\bf B258} (1985) 46--74.

\item{7.} P. Candelas, X.C. de la Ossa, P.S. Green and L. Parkes,
{\it A pair of Calabi--Yau manifolds as an exactly soluble superconformal
theory\/},  Nuc.\ Phys.\ {\bf B359} (1991) 21--74.

\item{8.} P. Candelas, X.C. de la Ossa, A.~Font, S.~Katz and D.R.~Morrison,
{\it Mirror symmetry for two-parameter models --- I\/}, to appear in
Nuc.\ Phys.\ B.

\item{} P. Candelas, A.~Font, S.~Katz and D.R.~Morrison,
{\it Mirror symmetry for two-parameter models --- II\/}, in preparation.

\item{9.} J. Carlson, M. Green, P. Griffiths and J. Harris,
{\it Infinitesimal variations of Hodge structure, I.\/}, Comp.\ Math.\
{\bf 50} (1983) 109--205.

\item{10.} H. Clemens,  {\it Homological equivalence, modulo algebraic
equivalence, is not finitely generated\/},  Publ.\ Math.\ IHES {\bf 58} (1983)
19--38.

\item{11.} H. Clemens,  {\it Some results on Abel-Jacobi mappings\/},
In: Topics in Transcendental Algebraic Geometry, Princeton Univ.\ Press,
Princeton, NJ 1984.

\item{12.} A. Collino, J. Murre, G. Welters, {\it On the family of conics
lying on a quartic threefold\/}, Rend.\ Sem.\ Mat.\ Univ.\ Pol.\ Torino
{\bf 38} (1980).

\item{13.} M. Dine, N. Seiberg, X.G. Wen, and E. Witten, {\it Non-perturbative
effects on the string world sheet\/},  Nucl.\ Phys.\ {\bf B278} (1987)
769--789.

\item{14.} G. Ellingsrud and S.A. Str\o mme,  {\it On the number of
twisted cubics on the general quintic threefold\/}, University
of Bergen preprint 63-7-2-1992.

\item{15.} G. Ellingsrud and S.A. Str\o mme, {\it The number of twisted cubic
curves on Calabi-Yau complete intersections\/},
Dyrkolbotn lectures, August 1993.

\item{16.} A. Font, {\it Periods and duality symmetries in Calabi-Yau
compactifications\/}, Nucl.\ Phys.\ {\bf B391} (1993) 358.

\item{17.} W. Fulton, {\it Intersection Theory\/},  Springer-Verlag,
Berlin-Heidelberg-New York 1984.

\item{18.} P. Griffiths and J. Harris, {\it Principles of Algebraic
Geometry\/}, Wiley-Interscience 1978.

\item{19.} R. Hartshorne, {\it Algebraic Geometry\/},  Springer-Verlag,
Berlin-Heidelberg-New York 1977.

\item{20.} S.~Hosono, A.~Klemm, S.~Theisen and S.-T.~Yau, {\it Mirror symmetry,
mirror map, and an application to Calabi-Yau hypersurfaces\/}, preprint,
HUTMP--93/0801.

\item{21.} T.~Johnsen and S.~Kleiman, {\it On Clemens' conjecture\/}, preprint.

\item{22.} S. Katz, {\it Degenerations of quintic threefolds and their
lines\/}, Duke Math.\ Jour.\ {\bf 50} (1983) 1127--1135.

\item{23.} S. Katz,  {\it On the finiteness of rational curves on quintic
threefolds\/},  Comp.\ Math.\ {\bf 60} (1986) 151--162.

\item{24.} S. Katz, {\it Lines on complete intersection threefolds with
$K = 0$\/},  Math.\ Zeit.\ {\bf 191} (1986) 293--296.

\item{25.} S. Katz,  {\it Iteration of multiple point formulas and applications
to conics\/},  Algebraic Geometry, Sundance 1986.  Lecture Notes in
Math.\ {\bf 1311} 147--155, Springer-Verlag, Berlin-Heidelberg-New York
1988.

\item{26.} S. Katz, {\it Curves on Calabi-Yau manifolds\/}, Dyrkolbotn
lectures, August 1993.

\item{27.} A.~Klemm and S.~Theisen,
{\it Considerations of one-modulus Calabi-Yau compactifications ---
 Picard-Fuchs equations, K\"ahler potentials and mirror maps\/},
Nucl.\ Phys.\ {\bf B389} (1993) 153.

\item{28.} K. Kodaira and D.C. Spencer,  {\it On deformations of complex
analytic structures, I.\/}, Annals of Mathematics {\bf 67} (1958) 328--401.

\item{29.} A.~Libgober, {\it Numerical characteristics of systems of
straight lines on complete intersections\/}, Math.\ Notes, {\bf 13} (1973)
51--56, Plenum Publishing, translated from Math.\ Zametki {\bf 13} (1973)
87--96.

\item{30.} A.~Libgober and J.~Teitelbaum, {\it Lines on Calabi-Yau complete
intersections, mirror symmetry, and Picard-Fuchs equations\/},
Int.\ Math.\ Res.\ Not.\ {\bf 1} (1993) 29--39.

\item{31.} D.R. Morrison, {\it Mirror symmetry and rational curves on
quintic threefold: a guide for mathematicians\/},  Jour.\ AMS {\bf 6}
(1993) 223--247.

\item{32.} D.R. Morrison, {\it Picard-Fuchs equations and mirror maps for
hypersurfaces\/}, these proceedings.

\item{33.} K. Oguiso, {\it Two remarks on Calabi-Yau Moishezon 3-folds\/},
preprint.

\item{34.} M. Reid, {\it Minimal models of canonical threefolds\/},  in:
Algebraic Varieties and Analytic Varieties.  Adv.\ Stud.\ Pure Math.\ {\bf 1}.
North-Holland, Amsterdam-New York and Kinokuniya, Tokyo, 1983,
131--180.

\item{35.} Y.~Ruan, {\it Topological sigma model and Donaldson type invariants
in Gromov theory\/}, pre- \linebreak print.

\item{36.} S.A. Str\o mme and D. van Straten, private communication
to the author by Str\o mme, 1990.

\item{37.} I. Vainsencher, {\it Elliptic quartic curves on a 5ic 3fold\/},
Proceedings of the 1989 Zeuthen Symposium, Cont.\ Math.\ {\bf 123} (1991)
247--257.

\item{38.} I. Vainsencher, {\it Enumeration of $n$-fold tangent hyperplanes
to a surface\/}, preprint.

\item{39.} X. Wu, {\it Chern classes and degenerations of hypersurfaces and
their lines\/},  Duke Math.\ J.\ {\bf 67} (1992) 633--652.

\bye